\documentclass[11pt,a4paper]{article}
\usepackage[T1]{fontenc}
\usepackage{authblk}
\usepackage{algorithm,algorithmic}
\usepackage{xcolor}
\usepackage{tikz}
\usepackage{cite}
\usepackage{comment}
\usetikzlibrary{positioning, arrows.meta, shapes.geometric}

\newcommand\keywords[1]{\textbf{Keywords}: #1}
\title{Chebyshev Spectral Neural Networks for Solving Partial Differential Equations}
\author[b]{Pengsong Yin}
\author[b]{Shuo Ling}
\author[a]{Wenjun Ying\thanks{Corresponding author.\\Email addresses:  wying@sjtu.edu.cn}}
\affil[a]{School of Mathematical Sciences, MOE-LSC and Institute of Natural Sciences,
Shanghai Jiao Tong University, Minhang,
Shanghai 200240, P.R. China.}
\affil[b]{School of Mathematical Sciences,
Shanghai Jiao Tong University,
Shanghai 200240, P.R. China}
 
\date{} 

\usepackage[left=1in, right=1in, top=0.5in, bottom=1.5in]{geometry}
\usepackage{indentfirst}  % 使用 indentfirst 宏包
\usepackage{graphicx} % Required for inserting images

\usepackage{amsmath,amsthm,amssymb,amsfonts}

\theoremstyle{definition}

\numberwithin{equation}{section}  % 公式编号包含章节信息

\date{\today}
\begin{document}
\maketitle
\begin{abstract}

The purpose of this study is to utilize the Chebyshev spectral method neural network (CSNN) model to solve differential equations. This approach employs a single-layer neural network wherein Chebyshev spectral methods are used to construct neurons satisfying boundary conditions. The study uses a feedforward neural network model and error backpropagation principles, utilizing automatic differentiation (AD) to compute the loss function. This method avoids the need to solve non-sparse linear systems, making it convenient for algorithm implementation and solving high-dimensional problems. The unique sampling method and neuron architecture significantly enhance the training efficiency and accuracy of the neural network. Furthermore, multiple networks enables the Chebyshev spectral method to handle equations on more complex domains. The numerical efficiency and accuracy of the CSNN model are investigated through testing on elliptic partial differential equations, and it is compared with the well-known Physics-Informed Neural Network(PINN) method.
\end{abstract}

% Keywords (see keywords.pdf file)
\keywords{Chebyshev spectral method neural network, differential equations, single-layer neural network, automatic differentiation, complex domains}

% Sections

\section{Introduction}

\indent

Partial differential equations (PDEs) play a crucial role in various fields of engineering and science. Many problems in mathematics, physics, engineering, economics, and other domains can be modeled using differential equations \cite{ricardo2020modern}. In most cases, obtaining analytical solutions for partial differential equations is challenging. Therefore, various numerical methods such as finite difference method (FDM) \cite{smith1985numerical}, finite element method (FEM) \cite{reddy1993introduction}, finite volume method (FVM) \cite{patankar2018numerical}, etc., are employed to solve these equations. While the aforementioned numerical methods provide good approximate solutions, they can be challenging for high-dimensional and complex domain problems.

Computational methods based on machine learning and neural networks have been widely researched in recent years. In fact, the increasing computational power and abundant data have made it possible to apply neural networks to increasingly complex problems. Furthermore, the versatility and simplicity of neural networks greatly solidify their powerful tool status in various domains. Meanwhile, the structure of neural networks is often reformulated to adapt to the diverse environmental requirements necessitated by different problems.

\subsection{Artificial Neural Networks for PDEs}
\indent

Recently, a class of Artificial Neural Networks(ANN) has been introduced, known as Physics-Informed Neural Networks (PINN) \cite{raissi2019physics, rudy2019data}. The training of PINN is achieved by minimizing the loss function, which, in addition to observational data, incorporates the governing equations, initial conditions, and boundary conditions. As is well-known, neural networks possess universal approximation capabilities \cite{ haykin1998neural, zurada1992introduction}. Therefore, following this process, approximate solutions to initial/boundary value problems may be effective.

Compared to traditional numerical methods, ANNs offer several advantages in providing approximate solutions. For example, such methods generally exhibit strong versatility, capable of solving both linear and nonlinear ordinary differential equations as well as partial differential equations. Additionally, many neural network approaches demonstrate certain superiority in solving high-dimensional and general domain problems.

While traditional neural networks have achieved remarkable results in solving differential equations, there are some drawbacks associated with using conventional neural network structures. Issues such as low solution accuracy and poor generality need to be overcome, and training may be challenging in certain situations, especially when the solution of the equation exhibits high-frequency characteristics \cite{wang2022and}. Therefore, many scientists are exploring approaches to address these challenges by incorporating mathematical prior knowledge, modifying network structures, training methods, and other techniques.

\subsection{Single Layer Artificial Neural Network Model}
\indent

Susmita Mall and S. Chakraverty proposed a fast learning single-layer neural network called Functional Link Artificial Neural Network (FLANN) to solve differential equations \cite{mall2014chebyshev, mall2017single}. In FLANN the hidden layer is replaced by a functional expansion block for enhancement of the input patterns using Chebyshev polynomials that is computationally more efficient than the multi-layer perceptron network \cite{patra2010chebyshev, zurada1992introduction}. While this single-layer Chebyshev network exhibits extremely fast optimization speed and, to some extent, has a mature approximation theory \cite{mason2002chebyshev}, using such a network structure alone for any PDE is not ideal. An important factor is that it does not take into account the need to satisfy the boundary conditions or initial conditions required for solving the equation.

To align the network structure more closely with the inherent properties of the equation, this paper proposes a new neural network solution for PDE based on the Chebyshev spectral method \cite{boyd2001chebyshev}. In the subsequent sections, we will refer to it as the Chebyshev spectral neural network(CSNN). Its advantage lies in its naturally satisfying the boundary conditions of the equation, as it avoids sampling at the boundaries, thereby reducing errors near the boundaries. Additionally, due to the excellent properties of spectral methods, such as exponential convergence for smooth functions, achieving high accuracy with fewer parameters of networks, and absence of phase errors \cite{shen2011spectral}, this method exhibits significant improvements in solution accuracy. In some cases, it even surpasses Physics-Informed Neural Networks (PINN) \cite{raissi2019physics}.

\subsection{Spectral Methods for Complex Domain Problems}
\indent

For the traditional spectral methods, addressing complex domain problems is crucial to expanding the application scope. Current mainstream solutions can be broadly categorized into two types: coordinate transformation and partitioning methods. Orszag introduced the concept of linear transformations \cite{orszag1979spectral}. Wang introduced general coordinates to spectral methods\cite{wang}, making it possible to compute for arbitrary shapes. Partitioning methods distribute the global nature of spectral methods to different regions to adapt to more complex geometries. Morchoishne first proposed spectral multi-domain methods \cite{morchoisne1982inhomogeneous}. Patera introduced spectral element methods \cite{patera1984spectral}, where spectral functions replace polynomials on each element, utilizing variational principles for finite element approximations. Although there are various methods available for solving complex domain problems, solving non-sparse matrices or intricate domain partitions remains a computationally expensive and challenging issue. CSNN employs variable substitution, transforming the polar coordinate set region that the equation needs to solve into a rectangular region for sampling and training. It uses an additional neural network to learn the transformed boundary values to achieve the solution, thus somewhat simplifying the computational complexity. While the form of the equation may undergo transformation, the powerful learning and generalization capabilities of neural networks enable effective learning even for nonlinear equations with low error. The subsequent chapters will provide a detailed explanation of this method.

The outline of this paper is as follows. Section \ref{Background Knowledge} introduces some prior knowledge about the Chebyshev spectral method. Section \ref{Chebyshev Spectral Method Neural Network} - \ref{solving Partial Differential Equations using CSNN} present the specific structure of the CSNN proposed by us and the algorithm for solving PDEs using it. Following that, section \ref{analysis} provides a certain error analysis. Section \ref{Numerical Experiments} is the numerical experiments section, validating the feasibility and efficiency of our method. The conclusion and discussion part is in section \ref{Conclusion and Discussion}.

\section{Preliminaries} \label{Background Knowledge}
This section introduces some knowledge about the traditional Chebyshev spectral method, serving as the foundation for CSNN.

\subsection{Chebyshev Polynomials}
\indent

The Chebyshev differential equation yields a set of orthogonal polynomials under the integral inner product on $[-1,  1]$ with the specific weight function $ w(x) = (1 - x^2)^{-\frac{1}{2}}$, which are known as the Chebyshev polynomials. The first two Chebyshev polynomials are known as
$$
T_0(x) = 1,  
$$
$$
T_1(x) = x.
$$

The well-known recursive formula can be employed to generate higher-order Chebyshev polynomials \cite{bhat2004numerical}
\begin{equation}
T_{n+1}(x) = 2xT_n(x) - T_{n-1}(x),    \qquad  n \geq 1,
\end{equation}
where $T_n(x)$ denotes n-th order Chebyshev polynomial and $ -1 < x < 1$ is the argument of the polynomials.

An important feature of the Chebyshev polynomials is that Chebyshev-Gauss-Lobatto quadrature points and wights can be expressed as follows
\begin{equation} \label{Chebyshev-Gauss-Lobatto quadrature points}
x_0 = 1, \quad x_N = -1, \quad x_j = \cos{\frac{\pi j}{N}},  \quad 1 \leq j \leq N-1,
\end{equation}
\begin{equation}
    \omega_0 = \omega_N = \frac{\pi}{2N}, \qquad \omega_j = \frac{\pi}{N}, \quad 1 \leq j \leq N-1.
\end{equation}

For Gauss quadrature defined with $\{ x_j, \omega_j\}_{j=0}^N$, define a discrete inner product in $C[a, b]$ and its associated norms:
\begin{equation} \label{discrete inner product}
    (u, v)_{N, \omega} = \sum_{j=0}^Nu(x_j)v(x_j)\omega_j,   \qquad \| u  \|_{N, \omega} = (u, u)_{N, \omega}^{\frac{1}{2}}
\end{equation}

\subsection{Basis Functions for Spectral Method} \label{Basis function for Spectral Method}
\indent

To introduce the construction method of basis functions, we exemplify with a (one-dimensional) two-point boundary-value problem. Without loss of generality, we consider the equation with homogeneous boundary conditions
\begin{equation}
\begin{aligned}
a_{-} u(-1)+b_{-} u^{\prime}(-1)&=0, \\
\qquad a_{+} u(1)+b_{+} u^{\prime}(1)&=0. 
\label{new boundary}
\end{aligned}
\end{equation}

For if the right-hand side is not zero, it can be always homogenized(see in Appendix A). Our objective is to identify a basis function of the following structure that satisfies \eqref{new boundary}
\begin{equation} \label{basis func}
\phi_k(x) = T_k(x) + a_kT_{k+1}(x) + b_{k}T_{k+2}(x).
\end{equation}
Since $T_k(\pm 1) = (\pm 1)^k$  and $T_k^{'}( \pm 1) = (\pm)^{k-1} k^2$, substitute \eqref{basis func} into  \eqref{new boundary} we have
$$
\left\{
\begin{aligned}
& \left(a_{+}+b_{+}(k+1)^2\right) a_k+\left(a_{+}+b_{+}(k+2)^2\right) b_k=-a_{+}-b_{+} k^2 \\
& -\left(a_{-}-b_{-}(k+1)^2\right) a_k+\left(a_{-}-b_{-}(k+2)^2\right) b_k=-a_{-}+b_{-} k^2
\end{aligned}
\right.
$$
Thus, we obtain the expressions for the computation of $a_k$ and $b_k$.
\begin{equation} \label{ak,bk}
\begin{aligned}
a_k= & -\left\{\left(a_{+}+b_{+}(k+2)^2\right)\left(-a_{-}+b_{-} k^2\right)\right. \\
& \left.-\left(a_{-}-b_{-}(k+2)^2\right)\left(-a_{+}-b_{+} k^2\right)\right\} / \text { DET }_k, \\
b_k= & \left\{\left(a_{+}+b_{+}(k+1)^2\right)\left(-a_{-}+b_{-} k^2\right)\right. \\
& \left.+\left(a_{-}-b_{-}(k+1)^2\right)\left(-a_{+}-b_{+} k^2\right)\right\} / \text { DET }_k,
\end{aligned}
\end{equation}
with
\begin{equation*} 
    \mathrm{DET}_k=2 a_{+} a_{-} - 2b_{+}b_{-}(k+1)^2(k+2)^2 + (a_{-}b_{+} - a_{+}b_{-})[(k+1)^2 + (k+2)^2].
\end{equation*}

Specifically, for the homogeneous Dirichlet conditions, substituting $a_k = 1$, $b_k = 0$, $c_k = -1$ into \eqref{ak,bk}, the basis functions can be easily obtained as
\begin{equation}
    \phi_i(x) = T_{i}(x) - T_{i+2}(x)
\end{equation}

For one-dimensional problems, we can directly take $\phi_k(x)$ as the basis for the approximation function space
\begin{equation} \label{1D Chebushev basis}
\mathcal{X}_N = span\{\phi_i(x):i = 0, 1, \dots, N-2\}.
\end{equation}

For a higher-dimensional cube (or hypercube) domain, we can compute polynomial basis functions for each coordinate corresponding to the homogeneous boundary conditions required on each dimension. Then, we can take the product of these basis functions to form new high-dimensional basis functions. For instance, an approximation function space for two-dimensional rectangle problem with homogeneous could be
\begin{equation} \label{basis}
\mathcal{X}_N = span\{\phi_i(x)\phi_j(y):i,j = 0, 1, \dots, N-2\}.
\end{equation}

\subsection{PDEs in Complex Geometries}
In this subsection, let us discuss how to solve equations in non-rectangular domains. Taking 2D problem as examples, assuming that we already have the polar coordinate representation for a complex domain in $\mathbb{R}^2$:
\begin{equation}
\begin{aligned}
& \mathrm{f}_1(\theta) \leq \mathrm{r} \leq \mathrm{f}_2(\theta), \\
& 0 \leq \theta<2 \pi,
\end{aligned}
\label{annular region}
\end{equation}
where $(r, \theta)$ are polar coordinates which can be transformed into the rectangle
\begin{equation} \label{rectangular}
\begin{aligned}
& -1 \leq z \leq 1, \\
& 0 \leq \theta<2 \pi,
\end{aligned}
\end{equation}
by the simple stretching transformation
\begin{equation}
z=2 \frac{r-f_1(\theta)}{f_2(\theta)-f_1(\theta)}-1 .
\end{equation}

The complexity of employing a coordinate transformation arises in the coefficients of the differential equation within the transformed domain\cite{orszag1979spectral}. Specifically, derivatives undergo transformation as per
$$
\begin{aligned}
& \left.\frac{\partial u}{\partial r}\right|_\theta=\left.\frac{2}{f_2(\theta)-f_1(\theta)} \frac{\partial u}{\partial z}\right|_\theta \\
& \left.\frac{\partial u}{\partial \theta}\right|_r=\left.\frac{\partial u}{\partial \theta}\right|_z-\left.\frac{\left[z\left(f_2^{\prime}-f_1^{\prime}\right)-\left(f_2^{\prime}+f_1^{\prime}\right)\right.}{f_2-f_1} \frac{\partial u}{\partial z}\right|_\theta .
\end{aligned}
$$

It is worth noting that nonlinear variable substitutions can sometimes introduce unexpected nonlinearities and singularities into the resulting differential equations. However, thanks to the powerful generalization and learning capabilities of neural networks, smaller errors can still be achieved when solving these differential equations. This will be discussed in detail in the Numerical Experimental sections.

\section{Chebyshev Spectral Method Neural Network(CSNN)} \label{Chebyshev Spectral Method Neural Network}

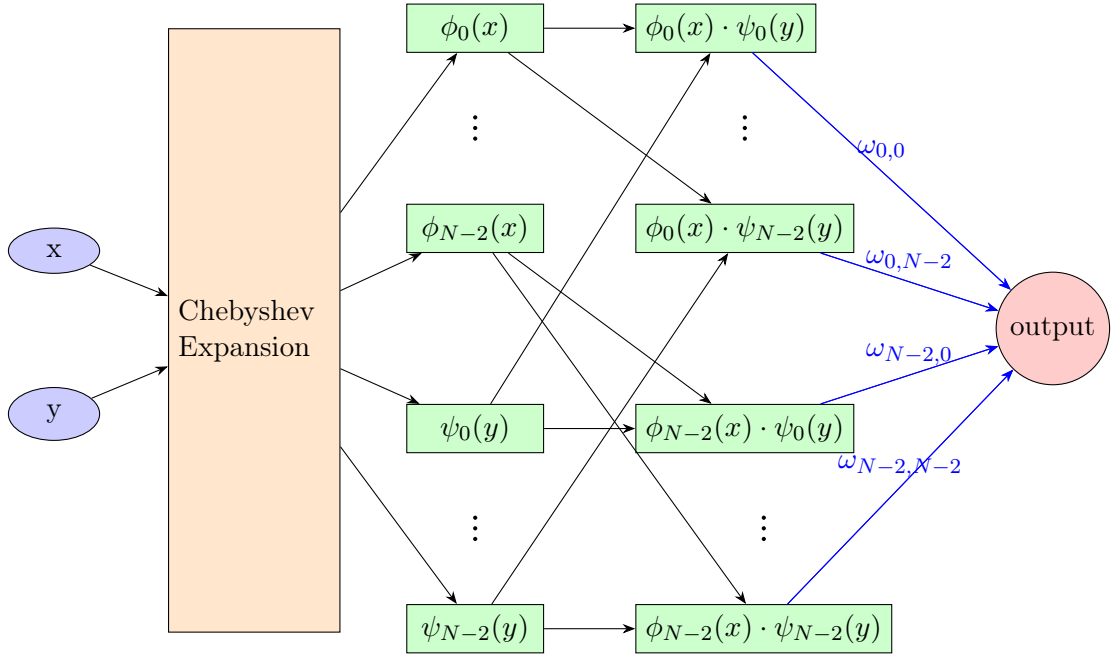
\begin{figure}[htbp]
    \centering
    \begin{tikzpicture}[
    % Define styles for different components
    input/.style={ellipse, minimum height=0.6cm, minimum width=1.2cm, draw=black, fill=blue!20},
    layer/.style={rectangle, minimum height=8cm, minimum width=2cm, draw=black, fill=orange!20},
    hidden/.style={rectangle, minimum height=0.6cm, minimum width=1.8cm, draw=black, fill=green!20},
    output/.style={circle, minimum height=1cm, minimum width=1cm, draw=black, fill=red!20},
    arrow/.style={->, >=Stealth},
    node distance=1.5cm and 1.5cm
    ]
    
    % Input layer
    \node[input]                  (input1) {x};
    \node[input, below=of input1] (input2) {y};
    
    % layer
    \path (input1) -- (input2) coordinate[midway] (midpoint);
    \node[layer, right=of midpoint, text width = 2cm] (layer1) {Chebyshev Expansion};
    
    % Hidden layer
    \node[hidden, right=2cm of layer1.north] (hidden1) {$\phi_0(x)$};
    \node[hidden, below=2cm of hidden1] (hidden2) {$\phi_{N-2}(x)$};
    \node[hidden, below=2cm of hidden2] (hidden3) {$\psi_0(y)$};
    \node[hidden, below=2cm of hidden3] (hidden4) {$\psi_{N-2}(y)$};
    \node[above=0.7cm of hidden2, font=\huge] {$\vdots$};
    \node[above=0.7cm of hidden4, font=\huge] {$\vdots$};
    
    \node[hidden, right=1.2cm of hidden1] (hidden21) {$\phi_0(x)\cdot\psi_0(y)$};
    \node[hidden, right=1.2cm of hidden2] (hidden22) {$\phi_0(x)\cdot\psi_{N-2}(y)$};
    \node[hidden, right=1.2cm of hidden3] (hidden23) {$\phi_{N-2}(x)\cdot\psi_{0}(y)$};
    \node[hidden, right=1.2cm of hidden4] (hidden24) {$\phi_{N-2}(x)\cdot\psi_{N-2}(y)$};
    \node[above=0.7cm of hidden22, font=\huge] {$\vdots$};
    \node[above=0.7cm of hidden24, font=\huge] {$\vdots$};
    
    % Output layer
    \path (hidden22) -- (hidden23) coordinate[midway] (hidden midpoint);
    \node[output, right=3.3cm of hidden midpoint] (output1) {output};
    
    % Connect layers
    \foreach \i in {1,2}
    \draw[arrow] (input\i) -- (layer1);
    \foreach \i in {1,2,3,4}
    \draw[arrow] (layer1) -- (hidden\i);
    \draw[arrow] (hidden1) -- (hidden21);
    \draw[arrow] (hidden3) -- (hidden21);
    \draw[arrow] (hidden1) -- (hidden22);
    \draw[arrow] (hidden4) -- (hidden22);
    \draw[arrow] (hidden2) -- (hidden23);
    \draw[arrow] (hidden3) -- (hidden23);
    \draw[arrow] (hidden2) -- (hidden24);
    \draw[arrow] (hidden4) -- (hidden24);

    \color{blue}
    \foreach \i in {1,2,3,4}
    \draw[arrow] (hidden2\i) -- (output1);
    \draw[arrow] (hidden21) -- node[pos=0.5, above] {$\omega_{0,0}$} (output1);
    \draw[arrow] (hidden22) -- node[pos=0.5, above] {$\omega_{0,N-2}$} (output1);
    \draw[arrow] (hidden23) -- node[pos=0.5, above] {$\omega_{N-2,0}$} (output1);
    \draw[arrow] (hidden24) -- node[pos=0.5, above] {$\omega_{N-2,N-2}$} (output1);

    \end{tikzpicture}
    
    \caption{The structure of CSNN model(only the blue lines contain the neural network weight parameters that need to be trained).}
    \label{CSNN}
\end{figure}

Figure \ref{CSNN} illustrates the structure of CSNN, which consists of an input layer with d inputs (input $\mathbf{x} \in \mathbb{R}^d$, in this figure and later discussion, we illustrate with the case where $d=2$), followed by neurons formed by basis functions based on Chebyshev polynomials, and a single output.

The CSNN model comprises three parts: The first part transforms the non-homogeneous problem into the corresponding homogeneous problem; the second part computes the Chebyshev basis functions that satisfy the conditionsas neurons; the third part discusses the learning and optimization algorithms of the neural network.

The first two parts have been elaborated in detail in Section \ref{Basis function for Spectral Method}. It should be noted that if the input elements satisfy the two-point boundary conditions(e.g., $z$ in \eqref{rectangular}), the basis functions can directly take the form of equation \eqref{basis func}. If the input elements are periodic elements introduced by variable substitution (e.g., $\theta$ in \eqref{rectangular}), due to their inherent periodicity, the spectral expansion is more suitable. We can choose the basis by $\phi_k(\theta) = \cos{k\theta} + \sin{k\theta}$, namely
\begin{equation} \label{trigonometric basis}
    \mathcal{\theta}_N = span\{\cos(i\theta): i = 0, 1, \dots, N-2\} \cup span\{\sin(i\theta): i = 0, 1, \dots, N-2\}.
\end{equation}
Experiments show that this choice indeed reduces the error more effectively than using Chebyshev basis functions for all cases directly.

The learning algorithm is employed to update the network parameters and minimize the loss function. In this case, the error backpropagation algorithm is utilized to update the weights of the CSNN. To accomplish this, the gradient of an loss function with respect to the network parameters $\mathbf{p}$ is computed. 

The network output with input data $\mathbf{x}$ and parameters (weights) $\mathbf{p}$ can be computed as:
$$
\mathcal{N}(\mathbf{x}; \mathbf{p})=\sum_{j=0}^{N-2}\sum_{i=0}^{N-2} w_{ij} \phi_{i}(x) \psi_{j}(y) .
$$
Where $\mathbf{x}=\left(x, y\right)$ represents the input data, $\phi_{i}(x)\psi_{j}(y)$ and $w_{ij}$ denote the basis function in \eqref{basis} and the weight vectors of CSNN, respectively.

% The weights of the CSNN can be adjusted using the principle of backpropagation:
% $$
% w_{ij}^{k+1}=w_{ij}^k+\Delta w_{ij}^k=w_{ij}^k+\left(-\eta \frac{\partial N(x, p)^k}{\partial w_{ij}^k}\right),
% $$
% where $\eta$ is the learning rate parameter, $k$ is the iteration step used for updating the weights as in conventional artificial neural networks, and $N(x, p)$ represents the error function.

Upon observation, it is evident that the function approximation space obtained using a single-layer CSNN neural network resembles that of spectral methods. In the subsequent numerical experiments, we also employ a single-layer network for experimentation. In fact, numerical results indicate that the single-layer CSNN exhibits high accuracy.

\section{Solving PDEs by CSNN} \label{solving Partial Differential Equations using CSNN}
The general form of representing ordinary differential equations (ODEs) or partial differential equations (PDEs) for the solution u($\mathbf{x}$) defined on a domain $\Omega \subset \mathbb{R}^d$ is as follows:
\begin{equation}
    f(\mathbf{x}, u(\mathbf{x}), \nabla u(\mathbf{x}), \nabla ^ 2 u(\mathbf{x}),  \dots, \nabla^{n} u(\mathbf{x})) = 0     \qquad \mathbf{x} \in \Omega
\end{equation}
where $f$ is a function that defines the structure of differential equation, $u(X)$ and $ \nabla $ denote the solution and differential operator respectively. With suitable boundary conditions
\begin{equation}
    \mathcal{B}(u, \mathbf{x}) = 0 \qquad  on \quad \partial \Omega
\end{equation}
where $\mathcal{B} (u, \mathbf{x})$ could be Dirichlet, Neumann, Robin boundary conditions. 

In the following discussion in section \ref{solving Partial Differential Equations using CSNN}, we will take \( d = 2 \) as an example while the case for higher dimensions is similar.

\subsection{PDEs in Rectangular Domains}
We begin by considering a simple rectangular domain. We first create a CSNN $\hat{u}(\mathbf{x}; \mathbf{\theta})$ as an approximation to the solution $u({\mathbf{x}})$. Then, we can take the derivatives of $\hat{u}(\mathbf{x}; \mathbf{\theta})$ with respect to the input $\mathbf{x}$ via automatic differentiation (AD) \cite{paszke2017automatic}. 

In fact, due to the specific structure of the CSNN, the network's output automatically satisfies the boundary conditions of the equation. Therefore, we only need to calculate the residual of the equation within the domain, without further considering the boundary conditions. This not only reduces the computational load of the loss function but also avoids the need to adjust the weights of the boundary and internal residuals as in traditional PINN methods.

As for sampling, if the training basis functions are Chebyshev polynomials\eqref{basis func}, then $N-2$ Chebyshev-Gauss-Lobatto quadrature points \eqref{Chebyshev-Gauss-Lobatto quadrature points} are taken in each dimension. If the training basis functions are Fourier polynomials, uniform sampling is directly used. This is a significant difference between CSNN and PINN. This sampling method greatly reduces the training workload.

We denote the sampling set as $\mathcal{T} = \{ \mathbf{x_1}, \mathbf{x_2},\dots, \mathbf{x_{\mathcal{|T|}}} \}$. To quantify the difference between the neural network $u({\mathbf{x}})$, we define the loss function as the weighted summation of the $L^2$ norm of residuals for the equation
\begin{equation}  \label{loss}
    \mathcal{L}(\mathbf{\theta}; \mathcal{T}) = \frac{1}{|\mathcal{T}|} \sum_{\mathbf{x} \in \mathcal{T}} \|  f(X, \hat{u}(X; \mathbf{\theta}), \nabla \hat{u}(X; \mathbf{\theta}), \nabla ^ 2 \hat{u}(X; \mathbf{\theta}),  \dots, \nabla^{n} \hat{u}(X; \mathbf{\theta})) \|_2^2
\end{equation}

Finally, by minimizing the loss function, we find the most suitable parameters $\mathbf{\theta}$. Considering that the loss with respect to $\mathbf{\theta}$ is highly nonlinear and non-convex \cite{blum1988training}, we typically utilize gradient-based optimizers such as gradient descent, Adam \cite{kingma2014adam}, and L-BFGS \cite{byrd1995limited} to minimize the loss function.

In conclusion, the process of solving PDEs on rectangular domains using the CSNN method is illustrated in Algorithm \ref{process1}.
\begin{algorithm}[!ht] \label{process1}
    \renewcommand{\algorithmicrequire}{\textbf{Input:}}
	\renewcommand{\algorithmicensure}{\textbf{Output:}}
	\caption{CSNN on rectangular domains}
    \label{power}
    \begin{algorithmic}[1] % 控制是否有序号
        \STATE For each dimension whose boundary condition is not of the type of periodic, homogenize
        boundary conditions into  \eqref{new boundary};
        \STATE  For each dimension whose boundary condition is not of the type of periodic, compute the values of $a_k, b_k, c_k$ using \eqref{ak,bk};
        \STATE Construct the Chebyshev spectral method neural Network using the computed $\phi_i$ \eqref{basis func};
        \STATE Compute the residual points $\mathcal{T}$ in each dimension;
        \STATE Train the neural network for a limited number of iterations with the loss function 
        $\mathcal{L}(\mathbf{\theta}; \mathcal{T})$ given by \eqref{loss}.
        % \STATE Estimate the mean $\mathrm{PDE}$ residual $\mathcal{E}_r$ by Monte Carlo integration, i.e., by the average of values at a set of randomly sampled dense locations $\mathcal{S}=\left\{\mathbf{x}_1, \mathbf{x}_2, \ldots, \mathbf{x}_{|\mathcal{S}|}\right\}:$
    \end{algorithmic}
\end{algorithm}

\subsection{PDEs in Complex Domains}
In this section, we will introduce how to solve PDEs on general parameterized domains, which is one of the major advantages of the CSNN method over traditional spectral methods and also one of our highlights.

\subsubsection{PDEs After Coordinate Transformation on Variables}

For more complex polar coordinate domains, when the inner boundary of the original region has prescribed values (such as in the case of an annular region \eqref{annular region}), the transformation poses no extra issues. This is because the boundary values on the rectangular region obtained from the transformation are determined by the inner and outer boundary conditions of the original region. Therefore, apart from the change in the form of the original equation due to the variable substitution and a new periodic boundary condition arise by $\theta$ direction, all other aspects are identical to solving PDEs on a rectangular domain(e.g. \eqref{rectangular}).

\begin{figure}[htbp]
	\centering
	\begin{minipage}{0.49\linewidth}
		\centering
		\includegraphics[width=0.6\linewidth]{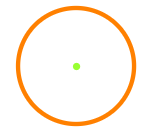}
		\caption{the original domain}
		\label{the complex domain}
	\end{minipage}
	\begin{minipage}{0.49\linewidth}
		\centering
		\includegraphics[width=0.6\linewidth]{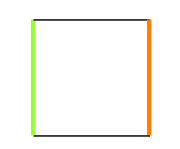}
		\caption{the parameter domain(with $r$ on the horizontal axis and $\theta$ on the vertical axis)}
		\label{the parameter domain}
	\end{minipage}
\end{figure}

However, when the original domain only has an outer boundary, such as
\begin{equation}
\begin{aligned}
& \mathrm{r} \leq \mathrm{f}_2(\theta), \\
& 0 \leq \theta<2 \pi.
\end{aligned}
\label{circular region}
\end{equation}
The values on one side of the transformed rectangular region are determined by the original region's boundary conditions, while the values on the other side become undetermined unknowns (Figure \ref{the complex domain} illustrates the original domain of the equation. Figure \ref{the parameter domain} illustrates the parameter domain in polar coordinates. The rectangular orange border corresponds to the boundary values of the original equation, while the green border corresponds to the unknown solution at the center of the original equation.)

% Traditional spectral methods may find it challenging to solve in such cases. But in CSNN, we can incorporate the boundary unknowns $u_0$ as parameters into the extra neural network, obtaining an approximate solution $\hat{u}(\mathbf{x}; \theta; u_0)$ to the equation. Alternatively, we can utilize a separate single-layer network $\hat{u_0}(\mathbf{x}; \theta)$ to learn the unknown values on the rectangular boundary. 

% In our experiments, we found that using the latter approach yields better results.  This is attributed to the ability to make correction to the $\hat{u_0}$ network independently during training, leading to substantial enhancements in computational accuracy and convergence speed.

Traditional spectral methods may find it challenging to solve in such cases. In CSNN, we can treat the boundary unknowns $u_0$ as an additional parameter $\hat{u_0}$ which can be trained during the training process together with $\hat{u}(\mathbf{x}; \theta)$. In addition, in the next subsection, we will introduce a correction technique for $u_0$ to assist in adjusting its value during the training process.

\subsubsection{Correction for $\hat{u_0}$ Network} \label{Correction}
Through multiple numerical experiments, we have found that after several iterations, the primary error in CSNN's solution to the differential equation arises from the discrepancy between the output values of the 
$\hat{u_0}$ network and the true solution at the origin. This error stems from never sampling the boundary points during training. We need to make certain correction to the 
$\hat{u_0}$ network to ensure that the equation is as close as possible to satisfying our original equation at the origin. 

Predefine an iteration threshold value $M$. When the number of iterations exceeds $M$, we perform a correction to $u$ every 100 iterations. Take 2D Possion equation as an example, we employ a five-point finite difference scheme for the correction
\begin{equation}
    -\frac{u(h, 0) - 2u(0, 0) + u(-h, 0)}{h^2} - \frac{u(0, h) - 2u(0, 0) + u(0, -h)}{h^2} \approx f(0, 0),
\end{equation}
where $h$ is a very small positive number. Thus, we obtain the correction form for the unknown value $u_0$ at the origin as
\begin{equation}
\bar{u_0} = \frac{\hat{u}(h, 0) + \hat{u}(-h, 0) + \hat{u}(0, h) + \hat{u}(0, -h)}{4} + h^2 f(0, 0).
\end{equation}
% where
% $$ 
% \left\{
% \begin{aligned}
% &\dot{u}(h, 0)  =  \hat{u}(2h - 1, -1) \\
% &\dot{u}(-h, 0)  =  \hat{u}(2h - 1, 0) \\
% &\dot{u}(0, h)  =  \hat{u}(2h - 1, -\frac{1}{2}) \\
% &\dot{u}(0, -h)  =  \hat{u}(2h - 1, \frac{1}{2}) \\
% \end{aligned}
% \right.
% $$
Here, $\hat{u}$ is the CSNN at the current stage, and $\bar{u}$ is the corrected prediction. In fact, we only need to update the value of $\bar{u}$ to the parameter of the single-parameter $\hat{u_0}$ network.

In conclusion, the process of solving PDEs on general parametric domains using the CSNN method is illustrated in Algorithm \ref{al2}.

\begin{algorithm}[!ht] \label{al2}
    \renewcommand{\algorithmicrequire}{\textbf{Input:}}
	\renewcommand{\algorithmicensure}{\textbf{Output:}}
	\caption{CSNN on general parametric domains}
    % \label{power}
    \begin{algorithmic}[1] % 控制是否有序号
        \STATE By variable substitution, the PDEs on general domains are transformed into new PDEs on rectangular domains;
        \STATE Determining the boundary conditions for the new PDEs and creating additional trainable parameters $\hat{u_0}$ caused by a new periodic boundary condition;
        \STATE Run Algorithm \ref{process1} on the new PDEs defined on the rectangular domain, and correct \( \hat{u_0} \) every certain number of training steps.
    \end{algorithmic}
\end{algorithm}

\section{Approximation Theory and Error Analysis for CSNN} \label{analysis}
Let $\mathcal{F}$ denote the family of all the functions that can be represented by the architecture of CSNN, 
 and u represent the true solution of the equation. We define $u_\mathcal{T} = arg \min_{f \in \mathcal{F}} \mathcal{L}(f; \mathcal{T})$ as the neural network whose loss is at global minimum on the trainning set $\mathcal{T}$. Since finding $u_\mathcal{T}$ through minimizing the loss is typically computationally challenging, we denote that our optimizer returns an approximate solution $\tilde {u_\mathcal{T}}$. Then, the total error $\mathcal{E}$ can be represented as \cite{bottou2007tradeoffs}
\begin{equation}
    \mathcal{E} := \| \tilde{u_\mathcal{T}} -  u \| \leq  \| \tilde{u_\mathcal{T}} - u_\mathcal{T}\| + \| u_\mathcal{T} - u\|.
\end{equation}

To simplify the problem, we assume that when the training samples are large enough, the global minimum point can be reached using the optimization algorithm of the neural network.i.e.
\begin{equation}
\| \tilde{u_\mathcal{T}} - u_\mathcal{T}\| \leq \epsilon
\end{equation}
 holds for any arbitrarily small $\epsilon$.

Here, for the approximation of $ \| u_\mathcal{T} - u\|$ we only consider the one-dimensional case. Suppose we need to solve the following model problem
\begin{equation} 
\begin{aligned}
    &-u^{''} + \alpha u = f,   \qquad in \quad I = (-1, 1)
    \label{model problem}  \\
    & u(1) = u(-1) = 0
\end{aligned}
\end{equation}
In this case, $\mathcal{F} = \{v \in P_N: v(-1) = v(1) = 0\}$. According to the assumptions, we have
\begin{equation} \label{one-dimensional equation}
\begin{aligned}
  &\alpha u_\mathcal{T}(x_j) - u_\mathcal{T}^{''}(x_j) = f(x_j),  \qquad  j = 1,2,\dots,N-1 \\
    &u_\mathcal{T}(x_0) = u_\mathcal{T}(x_N) = 0
\end{aligned}
\end{equation}
where $x_j = \cos(\frac{j \pi}{N})$. We can rewrite \eqref{one-dimensional equation} in a variational formulation with discrete inner product \eqref{discrete inner product} \cite{tang2006spectral}:
\begin{equation} \label{variational formulation}
    \alpha(u_\mathcal{T}, v_\mathcal{T})_{N, \omega} - (u_\mathcal{T}^{''}, v_\mathcal{T})_{N, \omega} = (I_{N}f, v_\mathcal{T})_{N, \omega}  \qquad \forall v_\mathcal{T} \in \mathcal{F}
\end{equation}
where $I_N:C(-1, 1) \to P_N$ is the interpolating operator associated with Gauss-Lobatto points.

According to  \cite{tang2006spectral}, for $u \in H_{\omega^{-\frac{3}{2},-\frac{3}{2}, *}}^m(I)$ and $f \in H_{\omega^{-\frac{1}{2},-\frac{1}{2}, *}}^k(I)$, we have
\begin{equation}
\begin{aligned}
& \left\|u-u_\mathcal{T}\right\|_{1, \omega} \lesssim N^{1-m}\left\|\partial_x^m u\right\|_{\omega^{m-3 / 2, m-3 / 2}}+N^{-k}\left\|\partial_x^k f\right\|_{\omega^{k-\frac{1}{2}, k-\frac{1}{2}}}, \quad m, k \geqslant 1, \\
& \left\|u-u_\mathcal{T}\right\|_\omega \lesssim N^{-m}\left\|\partial_x^m u\right\|_{\omega^{m-1 / 2, m-1 / 2}}+N^{-k}\left\|\partial_x^k f\right\|_{\omega^{k-\frac{1}{2}, k-\frac{1}{2}}}, \quad m, k \geqslant 1 .
\end{aligned}
\end{equation}
where $\omega$ represents the Chebyshev weight $(1 - x^2)^{-\frac{1}{2}}$, and the definition of the Jacobi-weighted Sobolev space can be found in Appendix B.

It indicates that if we ignore the error introduced by the neural network optimization algorithm, as long as the solution of the equation and the right-hand side $f$ are sufficiently smooth and the width of the neural network is wide enough, the approximate solution can converge to the true solution to any order.

\section{Numerical Experiments} \label{Numerical Experiments}
In this section, based on the aforementioned CSNN method, we validate the effectiveness and generality of the method by solving differential equations from one to four dimensions. All numerical experiments are conducted using the \textbf{NVIDIA GeForce RTX 3060(Laptop)} GPU, and the \textbf{Adam} optimizer is used to train the neural network. We begin with several simple examples.

We evaluate the accuracy by using the absolute $L^{\infty}$ error, absolute $L^2$ error and relative $L^2$ error defined as follows:
$$
\| e \|_{absolute-L^{\infty}} = \max_{i = 1, \cdots, K} | \hat{u}(x_i) - u^*(x_i) |,
$$
$$
\| e \|_{absolute-L^{2}} = \sqrt{\frac{\sum_{i = 1}^{K} | \hat{u}(x_i) - u^*(x_i) |^2}{K}} ,
$$
$$
\| e \|_{relative-L^{2}} = \frac{\sqrt{\sum_{i = 1}^{K} | \hat{u}(x_i) - u^*(x_i) |^2}}{\sqrt{\sum_{i = 1}^{K} |  u^*(x_i) |^2}} ,
$$

where $\hat{u}$ is the neural network prediction, $u^*$ is the true solution, and $K$ is the number of test points.

\subsection{The 1D Robin boundary condition differential equation}
Consider the following problem:
\begin{equation} \label{1D}
    u{''}(x) + xu^{'}(x) - u(x) = (24 + 5x)e^{5x} + (2+2x^2)\cos(x^2)-(4x^2+1)\sin(x^2)
\end{equation}
with the Robin boundary conditions:
\begin{equation}
    u(-1) - u^{'}(-1) = -4e^{-5} + \sin(1) + 2\cos(1), \qquad u(1) + u^{'}(1) = 6e^5 + \sin(1) + 2\cos(1).
\end{equation}

The exact solution is $u(x) = e^{5x}+\sin(x^2)$.

We take $N = 6, 8, 10, 12, 14$ to construct CSNN networks for computation. The epoch size for the training process of each network is set as $300 \times N$. The initial learning rate is set to $0.1$. A learning rate scheduler with a reduction factor of $0.7$ is utilized. This scheduler dynamically adjusts the learning rate if the loss fails to decrease for 600 consecutive iterations. 

The test points are uniformly sampled with 501 points in the interval from -1 to 1. The results are shown in Table \ref{tab:1d error}. 
\begin{table}[htbp]  \label{tab:1d error}
    \centering
    \begin{tabular}{|c|c|c|c|c|c|} \hline 
         $N$&  $6$&  $8$& $10$& $12$& $14$\\ \hline 
         Training Time (seconds) &10.44 &17.38 &23.51 &31.30 &41.21\\ \hline 
        Absolute $L^{\infty}$ Error  & 7.38E-2& 5.46E-3 &2.85E-4  &1.10E-5  &1.01E-6\\ \hline
    \end{tabular}
    \caption{Result for 1D Differential Equation}
    \label{Result for 1D Differential Equation}
\end{table}

It is observed that as the number of CSNN sampling points increases, the error of the solution decreases at a very high rate. In contrast, the slight increase in sampling points for PINNs has minimal impact on the equation's solution.

\subsection{The 2D Poisson equation on a rectangular domain}
In this part, we solve the two-dimensional Poisson equation with Dirichlet boundary conditions:
\begin{equation} \label{eq2DPisson_1}
 - \Delta u (x, y) =  f(x, y) \qquad (x, y) \in \Omega = [-1, 1]^2.
\end{equation}
The right-hand side term $f$ of the equation are determined by the true solution $u = e^{-x}\sin{(\pi y)}$. Let the boundary conditions of the equation be
\begin{equation}
\left\{\begin{array}{l}
\alpha_0(y)=u(-1, y), \\
\alpha_1(y)=u(1, y), \\
\beta_0(x)=u(x, -1), \\
\beta_1(x)=u(x, 1),
\end{array}\right.
\end{equation}
\begin{figure} 
\centering
\includegraphics[scale = 0.6]{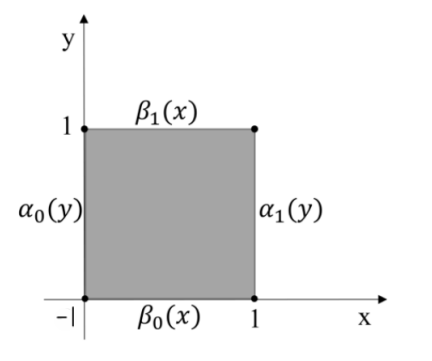}
\caption{the problem on the boundary}
\label{fig:axis}
\end{figure}
which Satisfy the continuity property
\begin{equation}
\left\{\begin{aligned}
\alpha_0(-1)&=u(-1,-1)=\beta_0(-1), \\
\alpha_0(1)&=u(1,1)=\beta_1(1), \\
\alpha_1(-1)&=u(1,-1)=\beta_0(1), \\
\alpha_1(1)&=u(1,1)=\beta_1(1) .
\end{aligned}\right.
\end{equation}

The trial solution is given by the expression,
\begin{equation}
    \psi(x, y)=A(x, y)+ \tilde{u}(x, y; \mathbf{\theta}).
\end{equation}
Here, $\tilde{u}(x, y; \mathbf{\theta})$ is the output obtained by the CSNN neural network, and $A(x, y)$ represents the construct function which satisfies the non-homogeneous boundary conditions, where
\begin{equation}
\begin{aligned}
    A(x, y) = &\left( \frac{1 - x}{2} \alpha_0(y) + \frac{1 + x}{2} \alpha_1(y) \right)\\
     &+ \left( \frac{1 - y}{2} \left( \beta_0(x) - \left( \frac{1 - x}{2} \beta_0(-1) + \frac{1 + x}{2} \beta_0(1) \right) \right) \right)\\
     &+ \left( \frac{1 + y}{2} \left( \beta_1(x) - \left( \frac{1 - x}{2} \beta_1(-1) + \frac{1 + x}{2} \beta_1(1) \right) \right) \right)
\end{aligned}
\end{equation}

The optimization algorithm chosen is the Adam, with a learning rate 0.1, and the number of iterations set to 2000. The test points are on a uniformly spaced grid of 300 cells on each dimension in the x and y space. The $L^{\infty}$ error and training duration for different network complexities N are shown in Table \ref{Result for 2D Poisson Equations}.
\begin{table}[htbp]
    \centering
    \begin{tabular}{|c|c|c|c|} \hline 
         $N$&  $6$&  $8$& $10$\\ \hline 
         Training Time (seconds)&  $19.16$&  26.18& 38.67\\ \hline 
         Absolute $L^\infty$ Error&  4.466E-4&  1.019E-5& 8.345E-7\\ \hline 
    \end{tabular}
    \caption{Result for 2D Poisson Equations}
    \label{Result for 2D Poisson Equations}
\end{table}

\begin{figure}[htbp]
  \centering
  \includegraphics[width=0.95\textwidth]{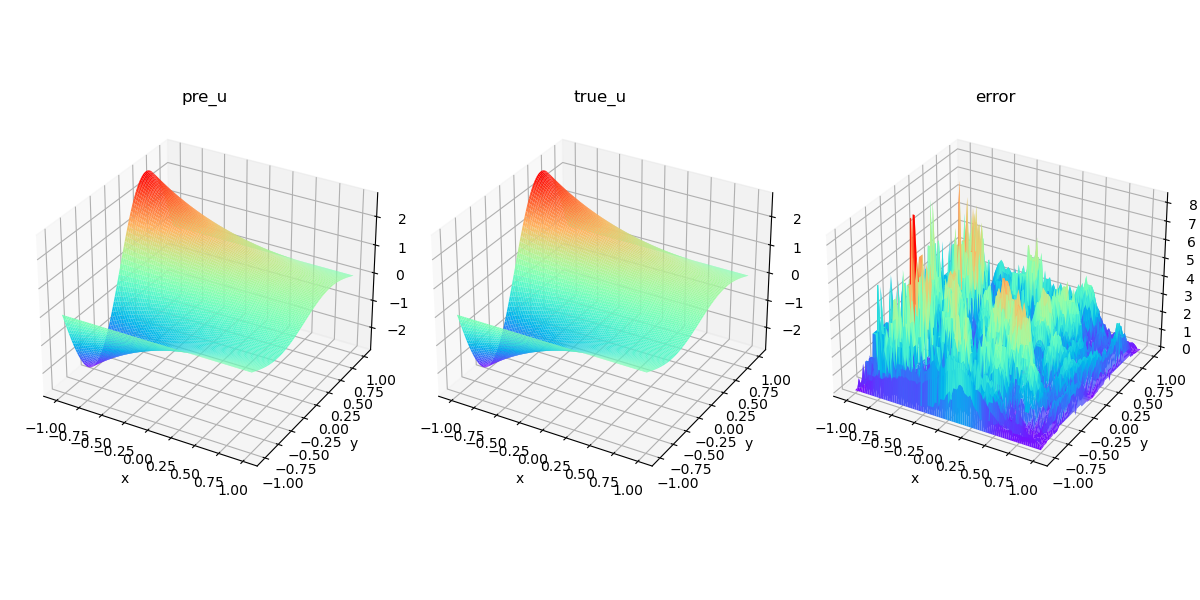}  
  \caption{predict for 2D Poisson equation(the left image shows the predicted values by the CSNN network, the middle image displays the accurate solution of the equation, and the right image represents the computed error).}
  \label{fig: Result for 2D Poisson Equations}
\end{figure}

A major advantage of the CSNN network is the selection of a large learning rate and the rapid convergence speed. The computational results for N = 10 are shown in Figure \ref{fig: Result for 2D Poisson Equations}. The error distribution is relatively uniform across the entire domain, with the error gradually decreasing near the boundary until it approaches zero.

\subsection{The 2D Poisson equation on an annular domain} \label{2D Poisson Equation on annular Domain}
In this part, we consider solving Poisson equation on the annular region $\Omega =\{(r, \theta) | 0.5 \leq r \leq 1$, $0 \leq \theta \leq 2\pi\}$. The original equation is given as follows:
\begin{equation} \label{eq2DPisson_2}
 - \Delta u (x, y) =  f(x, y) \qquad (x, y) \in \Omega.
\end{equation}
The right-hand side term $f$ and the Dirichlet boundary conditions on the inner and outer boundaries are both determined by the true solution $u = e^x\cos{(\pi y)}$.

Take $x = \frac{s + 3}{4}\cos((t+1)\pi)$, $y = \frac{s + 3}{4}\sin((t+1)\pi)$, the above equation \eqref{eq2DPisson_2} can be transformed into
\begin{equation}
    -\frac{16}{(s+3)} u_s - 16 u_{ss} - \frac{16}{{\left(\pi \cdot (s + 3)\right)}^2} u_{tt} = f(\frac{{3 + s}}{4}  \cos\left((t + 1) \pi\right), \frac{3+s}{4} \cos\left((t + 1)\pi\right)),
\end{equation}
where $t, s \in [ -1, 1]$.

In the s-direction, we employ Chebyshev basis functions (as in \eqref{1D Chebushev basis}), and the sampling points are still chosen as the Lobatto-Gauss-Legendre points \eqref{Chebyshev-Gauss-Lobatto quadrature points} (with N=10). For the t-direction, considering its periodicity, we use trigonometric basis functions (as in \eqref{trigonometric basis}), and the sampling points range uniformly from -1 to 1 (with 500 sampling points).

We denote the sampling set in the parameter space as $\mathcal{T} = \{ \mathbf{x_1}, \mathbf{x_2},\dots, \mathbf{x_{\mathcal{|T|}}} \}$. Constructing the loss function
\begin{equation}  
    \mathcal{L}(\mathbf{\theta}; \mathcal{T}) = \frac{1}{|\mathcal{T}|} \sqrt{ \sum_{\mathbf{x} \in \mathcal{T}} 
    | \frac{16}{(s+3)} u_s + 16 u_{ss} + \frac{16}{{\left(\pi \cdot (s + 3)\right)}^2} u_{tt} + f  |^2}
\end{equation}

\begin{figure}[htbp]
  \centering
  \includegraphics[width=\textwidth]{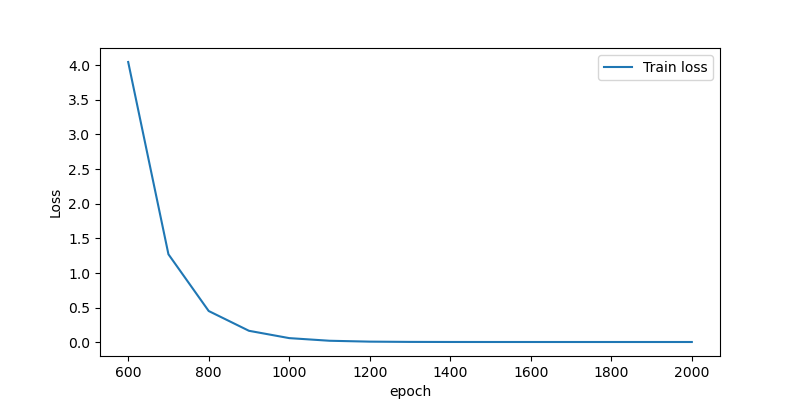}  
  \caption{loss for 2D Poisson Equations on annular Domain}
  \label{loss for 2D Poisson Equations}
\end{figure}

The initial learning rate is set to $0.01$, and a learning rate scheduler with a reduction factor of $0.6$ is utilized. This scheduler dynamically adjusts the learning rate if the loss fails to decrease for 500 consecutive iterations. As the training progresses, the loss function steadily decreases(see in Figure \ref{loss for 2D Poisson Equations}). It can be observed that even with a relatively small initial learning rate, the CSNN network converges quickly.

\begin{figure}[htbp]
  \centering
  \includegraphics[width=\textwidth]{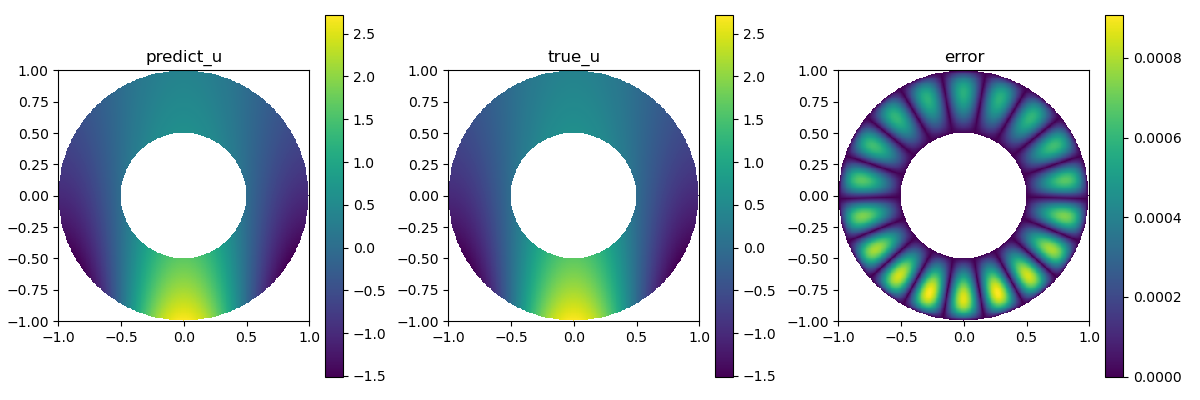}  
  \caption{The 2D Poisson equation on an annular domain(the left image shows the predicted values by the CSNN network, the middle image displays the accurate solution of the equation, and the right image represents the computed error.)}
  \label{fig: 2D Poisson Equations on annular Domain}
\end{figure}

The test points are on a uniformly spaced grid of 300 cells on each dimension in the s and t space. The solution of the equation is depicted in Figure \ref{fig: 2D Poisson Equations on annular Domain} and the absolute $L^{\infty}$ error reaches $9.07 \times 10^{-4}$ after 2000 training iterations. 

\subsection{The 2D Helmholtz equation on an annular domain} 
We consider the 2D Helmholtz equation on the $\Omega =\{(r, \theta) | 0.5 \leq r \leq 1$, $0 \leq \theta \leq 2\pi\}$
\begin{equation}
     \Delta u (x, y) + u(x, y) =  f(x, y).
\end{equation}
We test our method with $f$ and the Dirichlet boundary condition given by the exact solution $u(x, y) = e^x \sin(\pi y)$. Traditional methods for solving the 2D Helmholtz equation may encounter some difficulties. In finite difference, finite element, higher wave numbers may lead to numerical dissipation or numerical waveguide issues.

\begin{figure}[htbp]
  \centering
  \includegraphics[width=0.95\textwidth]{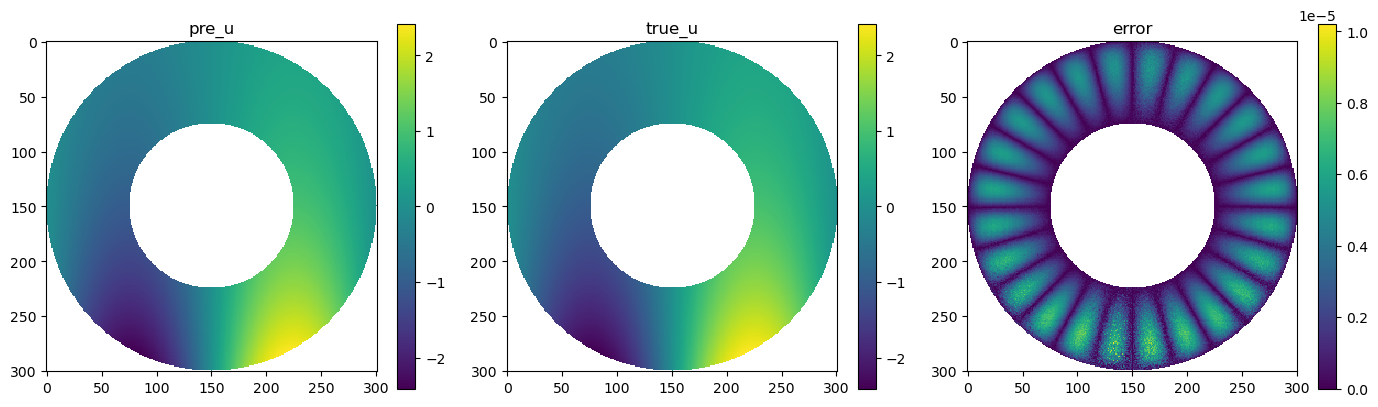}  
  \caption{The 2D Helmholtz equations on an annular domain(the left image shows the predicted values by the CSNN network, the middle image displays the accurate solution of the equation, and the right image represents the computed error.)}
  \label{fig: 2D MH Equations on annular Domain}
\end{figure}

In this example, we take $N = 13$, and the training parameters and test samples are the same as the previous subsection. The final numerical results are shown in Figure \ref{fig: 2D MH Equations on annular Domain} whose absolute $L_{\infty}$ error is $1.0192 \times 10^{-5}$. 

\begin{figure}[htbp]
  \centering
  \includegraphics[width=\textwidth]{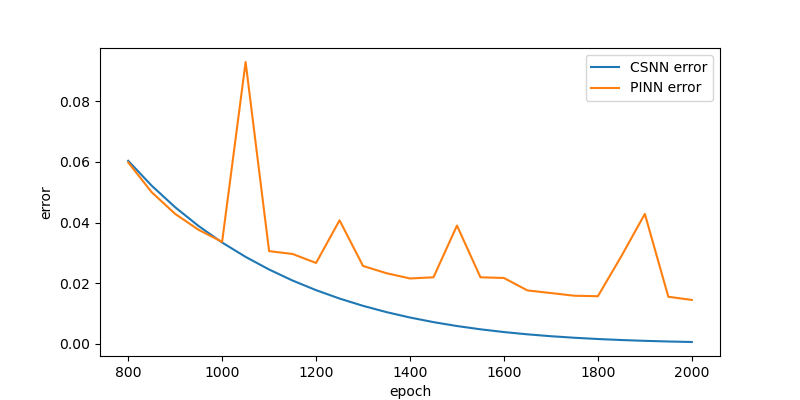}  
  \caption{Comparison with PINNs}
  \label{fig: compare}
\end{figure}

To compare the efficiency and accuracy differences between CSNN and PINNs, we applied the PINNs method to solve the examples in this section. The network of the PINNs method consists of 4 hidden layers, each with 128 neurons, using the $\tanh$ activation function. After 2000 iterations with the same learning rate $lr = 0.001$, the $L_{\infty}$ error of the solution given by the PINNs method is $0.014$, which is over 1000 times higher than CSNN. It can be observed that CSNN method has a higher convergence speed and can achieve higher accuracy as well. Simultaneously, during training, CSNN does not exhibit irregular oscillations in the solution error like PINNs; instead, the error steadily decreases (see in Fig \ref{fig: compare}).

\subsection{The 2D Poisson equation on a circular domain}
In this example, we consider solving the Poisson equation \eqref{eq2DPisson_2} on the circular region $\Omega =\{(r, \theta) | r \leq 1$, $0 \leq \theta \leq 2\pi\}$ with exact solution $u(x, y) = e^x \cos(\pi y)$. 

The initial learning rate is set to $lr=8 \times 10^{-3}$, and a learning rate scheduler with a reduction factor of $0.7$ is utilized. This scheduler dynamically adjusts the learning rate if the loss fails to decrease for 800 consecutive iterations. Additionally, we require an extra single-parameter neural network 
$\hat{u_0}$ to learn the value $u_0$ at the center of the circle. We set the initial learning rate for $\hat{u_0}$ as $6 \times 10^{-2}$ and decrease the step size by multiplying it by 0.8 if the loss fails to decrease for 600 consecutive iterations. 

\begin{figure}[htbp]
  \centering
  \includegraphics[width=\textwidth]{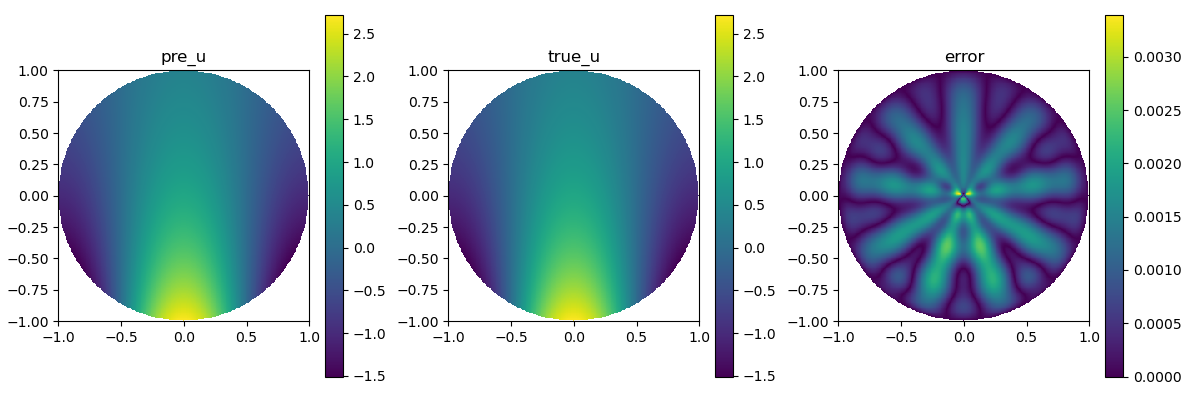}  
  \caption{the 2D Poisson equation on a circular domain(the left image shows the predicted values by the CSNN network, the middle image displays the accurate solution of the equation, and the right image represents the computed error.)}
  \label{fig: 2D MH Equations on circular domain}
\end{figure}

The test points are on a uniformly spaced grid of 300 cells on each dimension in the s and t space. After 10000 iterations of training, the $L^2$ error of the solution can reach $9 \times 10^{-4}$(see in Fig \ref{fig: 2D MH Equations on circular domain}).

However, it must be acknowledged that the feasibility of the CSNN method is strongly associated with the sizes of the learning rates for both networks. Incorrect step size settings can lead to slow convergence of the networks or even failure to reach the global minimum. To address this issue, we propose a network correction for $\hat{u}$ (as discussed in Section \ref{Correction}), where we apply adjustments to the network parameters to satisfy the equation conditions within a certain number of convergence steps. This measure significantly improves the resilience of the CSNN method to the learning rate settings.

\subsection{The 3D variable coefficients elliptic PDE} \label{3D Variable Coefficients Elliptic PDE}
\begin{figure}[htbp]
  \centering
  \includegraphics[width=0.65\textwidth]{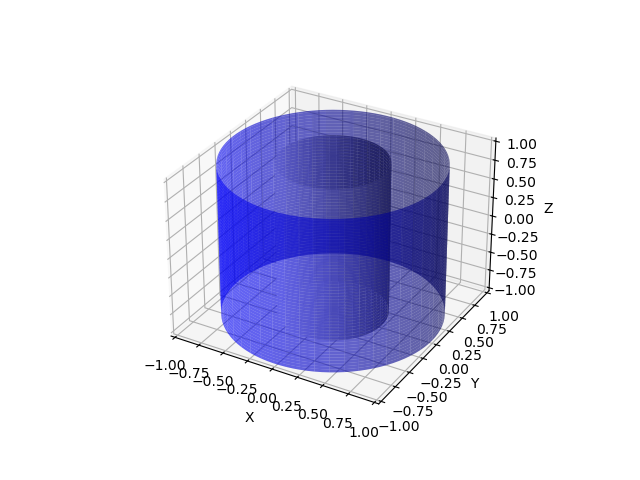}  
  \caption{The computational domain of Example \ref{3D Variable Coefficients Elliptic PDE}}
  \label{fig: 3D computational domain}
\end{figure}

For the three-dimensional case, we are interested in the numerical solution of the following form of the elliptic problem
\begin{equation} \label{3D PDE}
    \begin{aligned}
    &\nabla \cdot (\sigma(\mathbf{x})\nabla u) - \kappa(x) u = f(\mathbf{\mathbf{x}}) \qquad  & \text{in} \quad \Omega \\
    & u = g^{D}     &\text{on} \quad \partial \Omega
    \end{aligned}
\end{equation}
with $\sigma = \cos(\pi (x + y + z)) + 2$ and $\kappa = \sin(\pi (x + y + z)) + 2$.
The Dirichlet boundary data $g(x, y, z)$ and source data $f(x, y)$ are chosen such that the exact solution to the problem is given by $u(x, y, z) = e^{\frac{2}{3}x + \frac{2}{3}y + \frac{1}{3}z}\sin(\pi (x + y + z))$. The computational domain $\Omega$ is a rolled-up cylindrical shell with an inner radius of 0.5, an outer radius of 1, and a height ranging from -1 to 1(see Fig \ref{fig: 3D computational domain}).

If we take $x = r\cos \theta$, $y = r\sin \theta$, $z = z$, the aforementioned domain can be transformed into a standard rectangular region $\Omega^{'} = \{(r, \theta, z) \in [0.5, 1] \times [ 0, 2\pi) \times [-1, 1 ] \}$. And the differential operator is transformed into
\begin{equation}
\nabla \cdot\left(\sigma \nabla u(x, y, z)\right) = \frac{1}{r} \frac{\partial}{\partial r}(r \sigma \frac{\partial u}{\partial r}) + \frac{1}{r^2} \frac{\partial }{\partial \theta} (\sigma \frac{\partial u}{\partial \theta}) +  \frac{\partial }{\partial z} (\sigma \frac{\partial u}{\partial z})
\end{equation}

The initial learning rate is set to $0.01$, and a learning rate scheduler with a reduction factor of $0.7$ is utilized. This scheduler dynamically adjusts the learning rate if the loss fails to decrease for 100 consecutive iterations.  $r$ and $z$ directions are sampled at $N$ CGL points, and the $\theta$ direction is uniformly sampled with 100 points. The test points are on a uniformly spaced grid of 128 cells on each dimension in the $r, \theta, z$ space. In this experiment, $N = 8$, after training 78 seconds for 500 iterations, the $L_{\infty}$ error and $L_2$ error are $9.5367 \times 10^{-7}$ and $1.9601 \times 10^{-7}$, respectively. This demonstrates that even when solving nonlinear PDEs in irregular domains, CSNN still maintains a high level of accuracy.

\begin{figure}[htbp]
  \centering
  \includegraphics[width=\textwidth]{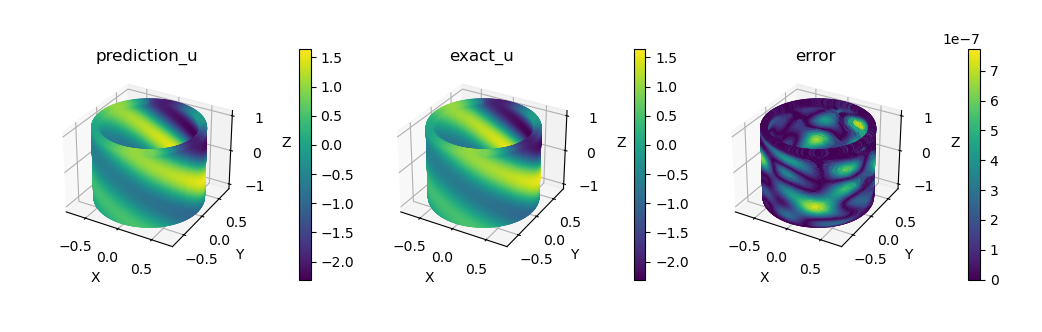}  
  \caption{The 3D variable coefficients elliptic(the left image shows the predicted values by the CSNN network, the middle image displays the accurate solution of the equation, and the right image represents the computed error.)}
  \label{fig: 3D Variable Coefficients Elliptic}
\end{figure}

As an example, for the cross-section at $r=0.7$, the numerical results and $L^{\infty}$ errors on the cylindrical surface are shown in Figure \ref{fig: 3D Variable Coefficients Elliptic}.

\subsection{The 4D Poisson Equation} \label{4D Poisson equations}
To illustrate the applicability of CSNN in higher dimensions, we solve the Poisson equation on a four-dimensional hypercube, testing the speed and accuracy of the method and comparing it with the traditional PINNs approach. The Poisson equation is given as following
\begin{equation} \label{eq4DPisson_1}
\begin{cases} - \Delta u (\mathbf{x}) = 4 \pi^2 sin(\pi x_1) sin(\pi x_2) sin(\pi x_3) sin(\pi x_4) & \mathbf{x} \in \Omega = [-1, 1]^4 \\ 
u(\mathbf{x}) = 0 & \mathbf{x} \in \partial \Omega \end{cases}.
\end{equation}
The exact solution is $u(\mathbf{x} = (x_1, x_2, x_3, x_4)) = sin(\pi x_1)sin(\pi x_2) sin(\pi x_3) sin(\pi x_4)$ which can be verified easily.

The result for our method on different values of $N$ after training 1000 epochs is given in Table \ref{Result for 4D Poisson Equations}. The test points are on a uniformly spaced grid of $11$ cells on each dimension in the 4D hypercube.
\begin{table}[htbp]
    \centering
    \begin{tabular}{|c|c|c|c|} \hline 
         $N$&  $6$&  $8$& $10$\\ \hline 
         Training Time (seconds)&  $79.20$&  90.71& 142.3\\ \hline 
         Absolute $L^2$ Error&  3.324E-2&  7.419E-4& 3.284E-5\\ \hline 
 Relative $L^2$ Error& 4.063E-2& 9.069E-4&4.014E-5\\ \hline
    \end{tabular}
    \caption{Result for 4D Poisson Equations}
    \label{Result for 4D Poisson Equations}
\end{table}
From the table, it can be observed that our method performs excellently in solving this PDE, especially for $N = 8$ and $10$. It demonstrates impressive speed and accuracy in numerical solutions for PDEs in high-dimensional spaces.

In addition, it is worth mentioning that the performance of the PINNs method in this case is catastrophic. Specifically, we employe a fully connected neural network with 3 hidden layers, each containing 128 neurons. The hyper-parameter $\beta$ for balancing PDE residual and boundary residual is set to 1, and the optimization scheme involves the widely used \textbf{Adam} optimizer. We randomly sample 65,536 points within $\Omega$ every 100 training steps to compute the corresponding $L^2$ norm inside the domain. Additionally, 1,024 points are randomly sampled along the boundaries to compute the corresponding $L^2$ norm on the boundary. Despite such substantial computational expenses, there is no corresponding improvement. Even after training for more than 10,000 iterations over 3,000 seconds, the loss function remains above 0.8. This results in the neural network failing to adequately fit the desired solution, with an absolute $L^2$ error reaching 0.88. This highlights the superiority of our method.

\section{Conclusion and Discussion} \label{Conclusion and Discussion} 
We propose a novel neural network framework for solving PDEs. Our experience with this method suggests that it offers the following advantages:
\begin{itemize}
  \item [1)] 
  High convergence rates and capability to fit high-frequency solutions.
  \item [2)]
  Reduced sampling point requirements, eliminating the need for boundary sampling points.
  \item [3)]
  Training using only a single-layer neural network greatly reduces the number of parameters, making convergence analysis feasible.
    \item [4)]
    Effective computation of mixed boundary conditions.
    \item [5)]
    Capable of solving equations in non-rectangular domains.
\end{itemize}
There are several disadvantages that need to be addressed in the future work:
\begin{itemize}
  \item [1)]  
In solving problems in complex domains, it is currently challenging to effectively determine the learning rates for each neural network.
  \item [2)]
  Inability to solve problems in non-parametrizable complex domains.
\end{itemize}

In future research, we will focus on solving more complex equations, such as the Stokes equations and the elasticity equations. We will also attempt to apply our methods to interface problems, hoping to achieve good numerical results. Additionally, more flexible and complex boundaries are also of interest to us. We hope to apply our method to any form of complex boundaries, not just parametric boundaries.

\section*{Acknowledgements} 
This work is financially supported by the Strategic Priority Research Program of Chinese Academy of Sciences(Grant No. XDA25010405). It is also partially supported by the National Key R \& D Program of China, Project Number 2020YFA0712000, the National Natural Science Foundation of China (Grant No. DMS-11771290) and the Science Challenge Project of China (Grant No. TZ2016002). 

\newpage
\section*{Appendix A} 
\renewcommand{\theequation}{\arabic{equation}} % 重定义公式编号格式，使其不包含章节号
Assume that U is the solution to a two-point boundary-value problem satisfying the following boundary conditions
\begin{equation} 
\begin{aligned}
a_{-} U(-1)+b_{-} U^{\prime}(-1) &= c_{-},  \\
\quad a_{+} U(1)+b_{+} U^{\prime}(1) &= c_{+}. 
\end{aligned}
\label{boundary-value}
\end{equation}

This includes in particular the Dirichlet $\left(a_{ \pm}=1\right.$ and $\left.b_{ \pm}=0\right)$, the Neumann $\left(a_{ \pm}=0\right.$ and $\left.b_{ \pm}=1\right)$, and the mixed $\left(a_{-}=b_{+}=0\right.$ or $a_{+}=b_{-}=0$ ) boundary conditions. The following explains how to homogenize it.

\begin{itemize}

\item Case $1 \quad a_{ \pm}=0$ and $b_{ \pm} \neq 0$. We set $\bar{u}=\beta x^2+\gamma x$, where $\beta$ and $\gamma$ are uniquely determined by asking $\tilde{u}$ to satisfy \eqref{boundary-value}
\begin{equation} \label{case1}
\begin{aligned}
-2 b_{-} \beta+b_{-} \gamma &=c_{-}, \\
\quad 2 b_{+} \beta+b_{+} \gamma&=c_{+}. 
\end{aligned}
\end{equation}

\item Case $2 \quad a_{-}^2+a_{+}^2 \neq 0$. We set $\tilde{u}=\beta x+\gamma$, where $\beta$ and $\gamma$ can again be uniquely determined by asking that $\tilde{u}$ to satisfy \eqref{boundary-value}
\begin{equation} \label{case2}
\begin{aligned}
\left(-a_{-}+b_{-}\right) \beta+a_{-} \gamma&=c_{-}, \\
\quad\left(a_{+}+b_{+}\right) \beta+a_{+} \gamma&=c_{+}.
\end{aligned}
\end{equation}

\end{itemize}

We now set $u=U-\tilde{u}$. Then $u$ satisfies the new equation with the homogeneous boundary conditions
\begin{equation}
\begin{aligned}
a_{-} u(-1)+b_{-} u^{\prime}(-1)&=0, \\
\qquad a_{+} u(1)+b_{+} u^{\prime}(1)&=0. 
\end{aligned}
\end{equation}

\section*{Appendix B} 
Let $I=(-1,1)$ and the Jacobi weight function of index $(\alpha, \beta)$ by $\omega^{\alpha, \beta}(x)=(1-x)^\alpha(1+x)^\beta$. We define the Jacobi weighted Sobolev spaces:
\begin{equation*}
\begin{aligned}
& L_{\omega^{\alpha, \beta}}^2(I)=\left\{u: \int_I u^2 \omega^{\alpha, \beta} \mathrm{d} x<+\infty\right\}, \\
& H_{\omega^{\alpha, \beta}}^l(I)=\left\{u \in L_{\omega^{\alpha, \beta}}^2(I): \partial_x u, \cdots, \partial_x^l u \in L_{\omega^{\alpha, \beta}}^2(I)\right\}, \\
& H_{\omega^{\alpha, \beta}, *}^m(I):=\left\{u: \partial_x^k u \in L_{\omega^{\alpha+k, \beta+k}}^2(I), \quad 0 \leqslant k \leqslant m\right\},
\end{aligned}
\end{equation*}

The norms in $L_{\omega^{\alpha, \beta}}^2(I)$ and $H_{\omega^{\alpha, \beta}}^l(I)$ will be denoted by $\|\cdot\|_{\omega^{\alpha, \beta}}$ and $\|\cdot\|_{l, {\omega^{\alpha, \beta}}}$, respectively. $H_{\omega^{\alpha, \beta}, *}^m(I)$ equipped with the inner product and norm
\begin{equation*}
(u, v)_{m, \omega^{\alpha, \beta}, *}=\sum_{k=0}^m\left(\partial_x^k u, \partial_x^k v\right)_{\omega^{\alpha+k, \beta+k}}, \qquad\|u\|_{m, \omega^{\alpha, \beta}, *}=(u, u)_{m, \omega^{\alpha, \beta}, *}^{\frac{1}{2}} .
\end{equation*}

% References
\newpage
\bibliographystyle{plain}
\bibliography{samples.bib}

\clearpage

\end{document}